\documentclass[a4paper, 10pt]{article}

\usepackage{geometry}
\usepackage{amsmath, amsthm, amssymb,amsfonts}
\usepackage[english, frenchb]{babel}
\usepackage{enumerate}
\usepackage{xcolor}
\usepackage[colorlinks=true, linktoc=page, citecolor=blue, linkcolor=blue, urlcolor=blue]{hyperref}

\usepackage{mathptmx}
\usepackage{setspace}
\numberwithin{equation}{section}

\usepackage{titlesec}
\renewcommand\thepart{\textbf{Part \Roman{part}.}}
\titleformat{\part}
  {\sc  \center}
  {\thepart}{1em}{}

\titleformat{\section}
  {\bf \normalsize \center}
  {\thesection.}{1em}{}

\titleformat{\subsection}
  {\bf \normalsize}
  {\thesubsection.}{1em}{}

\titleformat{\subsubsection}
  {\bf \normalsize}
  {\thesubsubsection.}{1em}{}  %

\newtheorem{theorem}{Theorem}[section]
\newtheorem{proposition}[theorem]{Proposition}
\newtheorem{cor}[theorem]{Corollary}
\newtheorem{lemma}[theorem]{Lemma}
\theoremstyle{definition}

\theoremstyle{remark}
\newtheorem{remark}[theorem]{Remark}

\numberwithin{equation}{section}

\newcommand{\C}{\mathbb{C}}
\newcommand{\g}{\mathfrak{g}}
\newcommand{\ka}{\mathfrak{k}}
\newcommand{\pe}{\mathfrak{p}}
\newcommand{\norme}[1]{\left\Vert #1\right\Vert}
\newcommand{\ind}{\mathrm{Ind}}
\newcommand{\indk}{\ind_{K_{\chi}}^K}
\newcommand{\ad}{\mathrm{Ad}}
\newcommand{\End}{\mathrm{End}}

\title{On the analogy between real reductive groups and Cartan motion groups. I: The Mackey-Higson bijection}
\author{Alexandre Afgoustidis\footnote{CNRS, Institut Élie Cartan de Lorraine, Nancy \& Metz, France.  \newline \emph{Mathematics Subject Classification (2020)}: 22E50, 22E46}}
\date{}

\begin{document}
\selectlanguage{english}
\maketitle

\begin{abstract}
George Mackey suggested in 1975 that there should be analogies between the irreducible unitary representations of a noncompact reductive Lie group $G$ and those of its Cartan motion group $G_0$ $-$ the semidirect product of a maximal compact subgroup of $G$ and a vector space. He conjectured the existence of a natural one-to-one correspondence between ``most'' irreducible (tempered) representations of $G$ and ``most'' irreducible (unitary) representations of $G_0$. We here describe a simple and natural bijection between the tempered duals of both groups, and an extension to a one-to-one correspondence between the admissible duals.
\end{abstract}

\section{Introduction}

\subsection{Contractions of Lie groups and a conjecture of Mackey}

\noindent When $G$ is a Lie group and $K$ is a closed subgroup, one can use the linear action of $K$ on the vector space $V =\mathrm{Lie}(G)/\mathrm{Lie}(K)$ to define a new Lie group: the semidirect product $G_0=K \ltimes V$, known as the \emph{contraction of $G$ with respect to $K$}. The notion first arose in mathematical physics (see Segal \cite{Segal}, \.{I}n\"{o}n\"{u} and Wigner \cite{InonuWignerContractions}, Saletan \cite{Saletan} and the lecture by Dyson \cite{Dyson}): the Poincar\'e group of special relativity admits as a contraction the Galilei group of classical inertial changes, and it is itself a contraction of the de Sitter group of general relativity. 

It is a classical problem to try to determine whether there is a relationship between the representation theories of $G$ and $G_0$: for instance, the unitary (irreducible) representations of the Poincar\'e group famously yield  spaces of quantum states for (elementary) particles \cite{WignerPoincare}, and it is quite natural to wonder about the existence of ``nonrelativistic'' analogues in the representation theory of the Galilei group \cite{InonuWignerGalilee}. In many cases of physical interest, the observation that the contraction $G_0$ can be seen as the special point of a one-parameter family $(G_t)_{t \in \mathbb{R}}$ of Lie groups (where the parameter $t$ usually has some physical significance, and $G_t$ is isomorphic with $G$ for all $t \neq 0$) has led to attempts to exhibit some representations of $G_0$ as ``limiting cases'' of representations of $G$; early examples include \cite{InonuWignerGalilee},  \cite{Hermann}, \cite{Mukunda},  \cite{MickelssonNiederle}.

For most Lie groups, unitary representations do not behave well under the contraction. The parameters needed to identify an irreducible representation of $G$ are often visibly different from those determining an irreducible representation of $G_0$, and some of the differences are important for physics $-$ a consequence of the bad behaviour in the case of the Poincar\'e group is that the notion of (rest) mass has different meanings in special and Galilean relativity. 
 
Now, suppose $G$ is a reductive Lie group and $K$ is a maximal compact subgroup. The contraction $G_0$ is known in that case as the \emph{Cartan motion group} of $G$, probably as a tribute to \'Elie Cartan's study of symmetric spaces: the group $G$ acts by isometries on the negatively curved $G/K$,  while $G_0$ acts by Euclidean rigid motions on the (flat) tangent space to $G/K$ at the identity coset. 

The algebraic structures of $G$ and $G_0$ and their representation theories are quite different. Denote the unitary dual of a Lie group $\Gamma$ by $\widehat{\Gamma}$; then George~Mackey's early work on semidirect products \cite{MackeyPNAS49, MackeyActa} describes $\widehat{G_0}$ in very simple and concrete terms, but describing $\widehat{G}$ remains to this day an extremely deep problem. A complete understanding of the \emph{tempered} dual $\widehat{G}^{\mathrm{temp}}$ (the support of the Plancherel measure) was attained in the early 1980s, crowning tremendous efforts of Harish-Chandra and others that had begun in 1945.
 
In 1975, however, Mackey noticed surprising analogies, for several simple examples of reductive groups $G$, between his accessible description of $\widehat{G_0}$ and Harish-Chandra's subtle parametrization of $\widehat{G}^{\mathrm{temp}}$. In the examples studied by Mackey, ``large'' subsets of $\widehat{G_0}$ and $\widehat{G}^{\mathrm{temp}}$ could be described using the same parameters, and the classical constructions for some of the corresponding representations were reminiscent of one another. Quantum-mechanical considerations led him to believe that there was more to this than chance, and he went on to conjecture that there should be a natural one-to-one correspondence between ``large'' subsets of $\widehat{G_0}$ and $\widehat{G}^{\mathrm{temp}}$  \cite{MackeyConjecture}: 
\vspace{0.15cm}

\begin{center}
\begin{minipage}{0.99\textwidth}
\emph{The physical interpretation suggests that there ought to exist a ``natural'' one-to-one correspondence between almost all the unitary representations of $G_0$ and almost all the unitary representations of $G$ $-$ in spite of the rather different algebraic structure of these groups.}
\end{minipage}
\end{center}
\vspace{0.15cm}

Mackey's idea seems to have been widely considered overenthusiastic at the time. It is true that in the years immediately following his proposal, studying the relationship between $G$-invariant harmonic analysis on $G/K$ and $G_0$-invariant analysis on $G_0/K$  became a flourishing subject (see e.g. \cite{HelgasonTangentSpace, Rouviere}), with beautiful ramifications for all Lie groups \cite{Duflo77, KashiwaraVergne}, and the close relationship between the two kinds of harmonic analysis does call to mind Mackey's observations on the principal series. But few of the authors who pursued this subject explicitly referred to Mackey's suggestion. Exceptions include Dooley and Rice \cite{DooleyRice}, who established in 1985 that the operators for (minimal) principal series representations of $G$  do weakly converge, as the contraction is performed, to operators for a generic representation of $G_0$; and Weinstein \cite{Weinstein}, who wrote down a Poisson correspondence relating coadjoint orbits of $G$ attached to (minimal) principal series representations and generic coadjoint orbits for $G_0$.

But as one moves away from the principal series of $G$ to the other parts of the tempered dual, it must have seemed difficult in the 1970s to imagine that anything general could be extracted from Mackey's suggestions. The geometrical construction of representations in the deeper parts of $\widehat{G}^{\mathrm{temp}}$, for instance that of the discrete series, was then a burning subject; in contrast, the most degenerate parts of $\widehat{G_0}$ look somewhat trivial. Mackey had of course noticed the absence, for these deeper parts, of any clear geometrical relationship:
\vspace{0.15cm}

\begin{center}
\begin{minipage}{0.99\textwidth}
\emph{Above all [the conjectured analogy] is a mere coincidence of parametrizations, with no evident relationship between the constructions of corresponding representations.}
\end{minipage}
\end{center}

\subsection{Connection with the Connes--Kasparov isomorphism; Higson's work}

\noindent In the late 1980s and in the 1990s,  Mackey's conjecture found an echo in the study of group $C^\star$-algebras. The \emph{Baum--Connes conjecture}, in its ``Lie group'' version due to Connes and Kasparov, describes the K-theory of the reduced $C^\star$-algebra $C_r^\star(G)$ of a connected Lie group $G$  in terms of the representation ring of a maximal compact subgroup and Dirac induction. For real reductive groups, the conjecture was first established by A.~Wassermann \cite{Wassermann}; a later proof would follow from V.~Lafforgue's work \cite{Lafforgue}. Viewing the operator K-theory of $C^\star_r(G)$ as a noncommutative-geometry version of the topological K-theory of the tempered dual $\widehat{G}^{\mathrm{temp}}$, Paul Baum, Alain Connes and Nigel Higson pointed out in the 1990s that the Connes--Kasparov isomorphism can be reinterpreted as a statement that $\widehat{G}^{\mathrm{temp}}$ and $\widehat{G_0}$ share algebraic-topological invariants: see \cite[\S 4]{BaumConnesHigson}. The Baum--Connes--Kasparov isomorphism can thus be viewed as a cohomological reflection of Mackey's conjectured analogy. 

Nigel Higson later remarked that Mackey had in mind a measure-theoretic correspondence (which may have been defined only for ``almost every'' representation), but that the interplay with the Connes--Kasparov isomorphism brought the topologies of $\widehat{G}^{\mathrm{temp}}$ and $\widehat{G_0}$ into the game, albeit at the level of cohomology. Because the relationship between the Fell topologies and the Plancherel measure is nontrivial, he noted \cite{Higson2008, HigsonComplex} that the simplest way to reconcile both points of view was to guess that there might exist a natural one-to-one correspondence between \emph{every irreducible tempered representation} of $G$ and \emph{every unitary irreducible representation} of $G_0$.

After the turn of the century, Higson took up the issue in more detail and examined the case of \emph{complex} semisimple groups, in which the tempered dual $\widehat{G}^{\mathrm{temp}}$ is completely described by the principal series. He constructed in 2008  a natural bijection between $\widehat{G}^{\mathrm{temp}}$ and $\widehat{G_0}$ in the complex case \cite{Higson2008}, and used it to analyze the structures of $C^\star_r(G)$ and $C^\star_r(G_0)$ $-$ discovering building blocks for the two $C^\star$-algebras that fit together rigidly in the deformation from $C^\star_r(G)$ to $C^\star_r(G_0)$.  With the help of tools crafted for the purpose in the early 1990s by Connes and himself (see \cite{ConnesHigsonDeformations} and \cite[\S II.10.$\delta$]{ConnesLivreENG}), he was led to a new proof of the Connes--Kasparov isomorphism for complex semisimple groups that is both natural from the point of view of representation theory and elementary from the point of view of K-theory.

\subsection{Contents of this paper} We here construct a simple and natural bijection between $\widehat{G}^{\mathrm{temp}}$ and $\widehat{G_0}$,  for any {real} reductive group $G$ (more precisely, whenever $G$ is the group of real points of a connected reductive algebraic group defined over $\mathbb{R}$. The construction may work for a wider class of reductive groups, as we will see). 

David Vogan's notion of lowest $K$-types will play a crucial part: we shall imitate Mackey's classical parametrization of $\widehat{G_0}$ (recalled in \S \ref{sec:section2}) and build a Mackey-like parametrization of $\widehat{G}^{\mathrm{temp}}$ using Vogan's work on tempered irreducible representations with real infinitesimal character (see \S \ref{subsec:vogan}). 

Our bijection will turn out (see \S \ref{subsec:ktypes}) to preserve lowest $K$-types and be compatible with a natural notion of variation of continuous parameters in the representation theories of both groups. As we shall see (in \S \ref{subsec:admissible}), the Langlands classification of admissible representations provides an immediate extension to a one-to-one correspondence between the admissible duals of both groups. 

This paper's main result, Theorem \ref{correspondance}, can of course be viewed a classification theorem for irreducible tempered representations. Its  overall form may seem particularly simple  in view of the equivalent classifications of Knapp--Zuckerman or Vogan--Zuckerman. One should not forget, though, that our proofs entirely rely on results extracted from these very deep classifications.

With Mackey's conjecture established at the level of parameters, we shall turn elsewhere to the deformation $(G_t)_{t\in \mathbb{R}}$ to try to better understand the reasons for, and implications of, the existence of the correspondence. A first companion paper \cite{AADuke} gives a geometrical realization for the Mackey--Higson bijection, describing a natural deformation (at the level of representation spaces) of every irreducible tempered representation $\pi$ of $G$ onto the corresponding representation of $G_0$. A second companion \cite{AAJFA} focuses on the topologies on $\widehat{G}^{\mathrm{temp}}$ and $\widehat{G_0}$ and the behavior of matrix coefficients, establishing rigidity properties of the deformation at the $C^\star$-algebraic level and deducing a new proof of the Connes--Kasparov isomorphism. 

Together, these companion papers lead to a somewhat different perspective on the bijection's definition, which we may summarize in elementary terms as follows (I thank a referee for this formulation). The irreducible tempered representations with a given lowest $K$-type can be parametrized by a real vector space (modulo a finite group). Moving along a ray in one of these vector spaces towards infinity, and inspecting the matrix coefficient functions attached to lowest $K$-type vectors in a representation $\pi$ of $G$, one can see that these functions actually converge, after a simple rescaling along  the deformation $(G_t)_{t\in \mathbb{R}}$, to certain matrix coefficients for a representation of $G_0$. This is the one corresponding to $\pi$. Unlike the short  algebraic construction of this paper, the more geometrical approach does not provide a proof for the bijection; but it is very relevant to the existing applications.

These further developments are hopefully clear signs that there may be much more to the Mackey--Higson bijection than the unexpected ``coincidence of parametrizations'' to be established here.

\subsection*{Acknowledgements}

Daniel Bennequin encouraged my interest for this subject and offered extremely kind and helpful advice. Several conversations with Michel Duflo, Nigel Higson and Mich\`ele Vergne were very useful and important to me. Fran\c cois Rouvi\`ere and David Vogan read the manuscript of my thesis \cite{AAthese} very carefully; the present version owes much to their suggestions and  corrections. Jeremy Daniel heard preliminary (and false) versions of many of my ideas. I am very much indebted to all of them.

During my time as a PhD student, in which many of the results of this paper were obtained, I benefited from excellent working conditions at Universit\'e Paris-7 and the Institut de Math\'ematiques de Jussieu--Paris Rive Gauche. This work was completed later, at Universit\'e Paris-Dauphine: I am very grateful to my colleagues at CEREMADE for their support.

\section{Representations of the Cartan motion group}\label{sec:section2}

\subsection{Notation}
\label{subsec:notations}

Let $G$ be the group of real points of a connected reductive algebraic group defined over $\mathbb{R}$ (see for instance \cite{HermitianFormsSMF}, \S 3), and let $\g$ be its Lie algebra. (All Lie algebras in this paper are real, unless they come with a subscript $\mathbb{C}$.) 

Fix a maximal compact subgroup $K$ of $G$, write $\theta$ for the Cartan involution of $G$ with fixed-point-set $K$, and form the corresponding Cartan decomposition $\g=\ka\oplus\pe$ of $G$. Here $\ka$ is the Lie algebra of $K$ and $\pe$ is a linear subspace of $\g$ which, although not a Lie subalgebra, is stable under the adjoint action $\ad: K \to \End(\g)$. 

Recall that the semidirect product $K \ltimes \pe$ is the group whose underlying set is the Cartesian product $K \times \pe$, equipped with the product law
\begin{equation*} (k_1, v_1) \cdot_0 (k_2, v_2) := (k_1k_2, v_1 + \mathrm{Ad}(k_1) v_2) \quad \quad {(k_{1}, k_2 \in K, \quad v_{1},v_{2}\in \pe).}\end{equation*}
We will denote that group by $G_0$ and call it the \emph{Cartan motion group of $G$}.

Writing $V^\star$ for the space of linear functionals on a vector space $V$ and viewing $\pe^\star$ as the space of linear functionals on $\g$ which vanish on $\ka$,  note that in the coadjoint action $\ad^\star: G \to \End(\g^\star)$, the subspace $\pe^\star \subset \g^\star$ is $\ad^\star(K)$-invariant.

We fix once and for all an abelian subalgebra $\mathfrak{a}$ of $\g$ that is contained in $\pe$ and is maximal among the abelian subalgebras of $\g$ contained in $\pe$. We write $A$ for $\exp_G(\mathfrak{a})$.

When $\Gamma$ is a Lie group, we will write $\widehat{\Gamma}$ for its unitary dual and $\widehat{\Gamma}^\mathrm{adm}$ for its admissible dual (see for instance \cite[\S 0.3]{Vogan81} or \cite[\S 3.3]{Wallach}). 
For the motion group $G_0$, which is amenable, the support of the Plancherel measure is all of  $\widehat{G_0}$. For our reductive group $G$, the support of the Plancherel measure is the \emph{tempered} dual $\widehat{G}^{\mathrm{temp}}$: the set of representations $\pi \in \widehat{G}$ whose $K$-finite matrix elements lie in $\mathbf{L}^{2+\varepsilon}(G/Z_G)$ for every $\varepsilon>0$. (See \cite{CowlingHaagerupHowe}; here $Z_G$ is the center of $G$.) If $G/Z_G$ is noncompact, then the trivial representation of $G$ is not tempered, and so $\widehat{G}^{\mathrm{temp}} \neq \widehat{G}$.

\subsection{Mackey parameters}
\label{subsec:data}

\noindent When $\chi$ is an element of $\pe^\star$, we write $K_\chi$ for its stabilizer in the coadjoint action of $K$ on $\pe$; this compact group is usually disconnected. In physics, it is known as the ``little group at $\chi$''.

We define a \emph{Mackey parameter} (or: \emph{Mackey datum}) to be any pair $\delta=(\chi, \mu)$ where $\chi$ is an element of $\pe^\star$ and $\mu$ is an irreducible representation~of~${K_\chi}$.

We say that two Mackey parameters $(\chi, \mu)$ and $(\chi', \mu')$ are \emph{equivalent} when there exists an element $k$ in $K$ such that 
\begin{itemize}
\item[$\bullet$] $\chi=\ad^\star(k)\cdot \chi'$, and
\item[$\bullet$] $\mu$ is equivalent, as an irreducible representation of $K_\chi$, with  the representation $u \to \mu'(k^{-1}uk)$.
\end{itemize}
This defines an equivalence relation on the set of Mackey data; we will write $\mathcal{D}$ for the set of equivalence classes.

It may be useful to mention two elementary facts from structure theory (see \cite{KnappPrincetonBook}, sections V.2 and V.3): 

\begin{itemize}
\item[(i)] Every Mackey datum is equivalent to a Mackey datum $(\chi, \mu)$ in which $\chi$ lies in $\mathfrak{a}^\star$; furthermore, if $\chi$ and $\chi'$ are two elements of $\mathfrak{a}^\star$, and if $\mu$ and $\mu'$ are irreducible representations of $K_{\chi}$ and $K_{\chi'}$, then the Mackey data $(\chi, \mu)$ and $(\chi', \mu')$ are equivalent if and only if there is an element of the Weyl group $W = W(\g, \mathfrak{a})$ which sends $\chi$ to $\chi'$ and $\mu$ to $\mu'$.
\item[(ii)] The dimension of the $\ad^\star(K)$-orbit of $\chi$ in $\pe^\star$ is maximal (among all possible dimensions of $\ad^\star(K)$-orbits in $\pe^\star$) if and only if $\chi$ is \emph{regular}. All regular elements in $\mathfrak{a}^\star$ have the same $\ad^\star(K)$-stabilizer, \emph{viz.} the centralizer $M=Z_K(\mathfrak{a})$ of $\mathfrak{a}$ in $K$.
\end{itemize}


\subsection{Unitary dual of the motion group}
\label{subsec:dualG0}

\noindent We recall the usual way to build a unitary representation of $G_0$ from a Mackey parameter $\delta=(\chi, \mu)$. 

Consider the centralizer $L_0^\chi$ of $\chi$ in $G_0$ (for the coadjoint action). It is $K_\chi \ltimes \pe$, where $K_\chi$ is the little group of \S \ref{subsec:data}.

Out of $\delta$ and $\chi$, build an irreducible representation of this centralizer: fix a (finite-dimensional) carrier space $V$ for $\mu$, and define an action of $L_\chi^0$ on $V$, where $(\kappa, v) \in K_\chi \ltimes \pe$ acts through $e^{i\langle \chi, v\rangle} \mu(k)$. Write $\sigma=\mu \otimes e^{i\chi}$ for this representation of $L_\chi^0$. 

Now, observe the induced representation
\begin{equation} \mathbf{M}_0(\delta)=\ind_{L^0_\chi}^{G_0}(\sigma) = \ind_{K_\chi \ltimes \pe}^{G_0}\left(\mu\otimes e^{i\chi}\right).\end{equation}

 \begin{theorem}[Mackey \cite{MackeyPNAS49}] \leavevmode
 \begin{enumerate}[(a)]
 \item For every Mackey parameter $\delta$, the representation $\mathbf{M}_0(\delta)$ is irreducible.
 \item Suppose $\delta$ and $\delta'$ are two Mackey parameters. Then the representations $\mathbf{M}_0(\delta)$ and $\mathbf{M}_0(\delta')$ are unitarily equivalent if and only if $\delta$ and $\delta'$ are equivalent as Mackey parameters.
 \item By associating to any Mackey parameter $\delta$ the equivalence class of $\mathbf{M}_0(\delta)$ in $\widehat{G_0}$, one obtains a bijection between $\mathcal{D}$ and $\widehat{G_0}$. \end{enumerate}\end{theorem}

\begin{remark} \label{copiekhat} If $\chi$ is zero, then $K_\chi$ equals $K$ and the representation of $\mathbf{M}_0(\chi, \mu)$ of $G_0$ is the trivial extension of $\mu$ where $\pe$ acts by the identity. There is thus in $\widehat{G_0}$ a family of finite-dimensional representations; every member of the family is equal to $\mathbf{M}_0(0, \mu)$ for some $\mu$ in $\widehat{K}$.\end{remark}

\begin{remark} \label{constructionG0} In \S \ref{subsec:ktypes}, we will need information about the $K$-module structure of $\mathbf{M}_0(\chi, \mu)$. Recall that $\mathbf{M}_0(\chi, \mu)$ can be realized by fixing a $\mu(K)$-invariant inner product on $V$, then equipping the Hilbert space
\begin{equation*} \label{espaceG0} \mathbf{H}= \left\{ f \in \mathbf{L}^2(K, V)\ / \ \forall m \in K_{\chi}, \forall k \in K, \ f(km)=\mu(m^{- 1})f(k) \right\}\end{equation*}
 with the $G_{0}$-action in which 
\begin{equation*} \label{actionG0} g_{0} = (k,v) \text{ acts through } \pi_0(k,v): f \mapsto \left[ u \mapsto e^{i \langle Ad^\star(u)\chi, v\rangle} f(k^{- 1} u)\right];\end{equation*}
the restriction to $K$ of that action~is~a~classical~picture~for~$\ind_{K_\chi}^K(\mu)$.\end{remark}

\subsection{Admissible dual of the motion group}\label{subsec:admG0}

\noindent We now outline Champetier and Delorme's description of the admissible dual of $\widehat{G_0}^{\mathrm{adm}}$ (see \cite{ChampetierDelorme}; see also Rader \cite{Rader}). Recall that the action $\ad^\star$ of $K$ on $\pe^\star$ induces an action on the complexified space $\pe_\C^\star$, in which $K$ acts separately on the real and imaginary parts. We will still write $\ad^\star$ for it.

We define an \emph{admissible Mackey parameter} to be any pair $(\chi, \mu)$ in which $\chi$ lies in $\mathfrak{a}_\C^\star$ and $\mu$ is an irreducible representation of $K_\chi$. In contrast to the previous section, we require from the outset that $\chi$ lie in $\mathfrak{a}_\C^\star$: it is no longer true that $\ad^\star(K) \cdot \mathfrak{a}_\C^\star$ coincides with $\pe_\C^\star$.

We define on the set of admissible Mackey parameters the same equivalence relation as described in \S \ref{subsec:data} and write $\mathcal{D}^\mathrm{adm}$ for the set of equivalence classes of admissible Mackey data.

When $\delta=(\chi, \mu)$ is an admissible Mackey datum, we write $\mathbf{M}_0(\delta)$ for the representation  $\mathrm{Ind}_{K_{\chi} \ltimes \pe}^{G_{0}}(\mu \otimes e^{\chi})$ of $G_{0}$. This is an admissible representation of $G_0$ (for details on admissible representations of $G_0$, and equivalences between them, see \cite{ChampetierDelorme}). This representation is unitary if and only if the real part of $\chi$ is zero; in that case it coincides with the representation defined in \S \ref{subsec:dualG0} using the imaginary part of $\chi$.

\begin{theorem}[\cite{ChampetierDelorme}, Th\'eor\`eme A] \leavevmode
\begin{enumerate}[(a)]
\item For every  admissible Mackey parameter $\delta$, the representation $\mathbf{M}_0(\delta)$ is irreducible.
\item Suppose $\delta$ and $\delta'$ are two admissible Mackey parameters. Then the admissible representations $\mathbf{M}_0(\delta)$ and $\mathbf{M}_0(\delta')$ are equivalent if and only if $\delta$ and $\delta'$ are equivalent as admissible Mackey parameters.
 \item By associating to any admissible Mackey parameter $\delta$ the equivalence class of $\mathbf{M}_0(\delta)$ in $\widehat{G}^{\mathrm{adm}}$, one obtains a bijection between $\mathcal{D}^{\mathrm{adm}}$ and $\widehat{G_0}^{\mathrm{adm}}$. 
 \end{enumerate}\end{theorem}

\section{The Mackey--Higson bijection for tempered representations}

\noindent In this section, we introduce a natural one-to-one correspondence between $\widehat{G_0}$ and $\widehat{G}^{\mathrm{temp}}$. 

The key step is to determine the subset of $\widehat{G}^{\mathrm{temp}}$ that should correspond to the ``discrete part'' of $\widehat{G_0}$ encountered in Remark \ref{copiekhat}, which is a copy of $\widehat{K}$ in $\widehat{G_0}$.

For some reductive groups, there is a natural candidate for that subset of $\widehat{G}^{\mathrm{temp}}$: the discrete series. If $G$ is a reductive group with nonempty discrete series, it is well-known from the work of Harish-Chandra, Blattner, Hecht--Schmid (see \cite{HechtSchmid}) that every discrete series representation $\pi$ comes with a distinguished element of $\widehat{K}$, the minimal  $K$-type of $\pi$ (see for example \cite[\S 1.5]{DufloSemBou}); furthermore, inequivalent discrete series representations of $G$ have inequivalent minimal $K$-types.

There are at least two obvious reasons not to use the discrete series for our construction, however: 
\begin{enumerate}[(i)]
\item there are reductive groups with no discrete series representations;
\item even when $G$ has a nonempty discrete series, there are elements of $\widehat{K}$ which cannot be realized as the minimal $K$-type of any discrete series representation.
\end{enumerate}
\indent To define a natural candidate for the class of representations to be associated in $\widehat{G}^{\mathrm{temp}}$ to the copy of $\widehat{K}$ in $\widehat{G_0}$, we turn to David Vogan's work on lowest $K$-types for other representations.

\subsection{Real-infinitesimal-character representations and a theorem of Vogan} \label{subsec:vogan}

\noindent Fix a reductive group $G$ and a maximal compact subgroup $K$. Given a Cartan subalgebra $\mathfrak{t}$ of $\ka$ and and a choice of positive root system $\Delta^+_c$ for the pair $(\ka_\C, \mathfrak{t}_\C)$, Vogan defines (\cite{Vogan77}) a positive-valued function $\mu \to \norme{\mu}_{\widehat{K}}$ on $\widehat{K}$: starting from an irreducible $\mu$ in $\widehat{K}$, one can consider its $\Delta^+_c$-highest weight $\tilde{\mu}$ and the half-sum $\rho_c$ of positive roots in $\Delta_c^+$; these are two elements of  $\mathfrak{t}_{\C}^\star$, and that space comes with a Euclidean norm $\norme{\cdot}_{\mathfrak{t}_{\C}^\star}$; one then defines $\norme{\mu}_{\widehat{K}} = \norme{\tilde{\mu}+2\rho_{c}}_{\mathfrak{t}_{\C}^\star}$.

When $\pi$ is an admissible (but not necessarily irreducible) representation of $G$, the \emph{lowest $K$-types of $\pi$} are the elements of $\widehat{K}$ that appear in the restriction of $\pi$ to $K$ and have the smallest possible norm among the elements of $\widehat{K}$ that appear in $\pi_{|K}$. Every admissible representation has a finite number of lowest $K$-types, and these do not depend on the choice of $T$ and $\Delta_c^+$. 

Given an irreducible tempered representation $\pi$ of $G$, we say $\pi$ has \emph{real infinitesimal character} when there exists a cuspidal parabolic subgroup $P$ of $G$, with Langlands decomposition $P=MAN$, and a discrete series representation $\sigma$ of $M$, so that $\pi$ is equivalent with one of the irreducible factors of $\ind_{P}^G(\sigma \otimes \mathbf{1})$. We will not need to explain the link with the notion of infinitesimal character and refer the reader to \cite{VoganLanglands}.

We write $\widehat{G}^{\mathrm{RIC}}$ for the subset of $\widehat{G}^{\mathrm{temp}}$ whose elements are the equivalence classes of irreducible tempered representations of $G$ with real infinitesimal character. I thank Michel Duflo for calling my attention to  that subset.

Vogan discovered (\cite{Vogan81}, see also \cite{Vogan79}) that if $\pi$ is an irreducible tempered representations of $G$ with real infinitesimal character, then $\pi$ has a \emph{unique} lowest $K$-type. He went on to prove that every $K$-type is the lowest $K$-type of exactly one representation in $\widehat{G}^{\mathrm{RIC}}$:

\begin{theorem}[Vogan \cite{Vogan81}; see \cite{VoganBranching}, Theorem 1.2]  \label{thmktypes} The map from $\widehat{G}^{\mathrm{RIC}}$ to $\widehat{K}$ which, to an irreducible tempered representation $\pi$ of $G$ with real infinitesimal character, associates the lowest $K$-type of $\pi$, is a bijection.
\end{theorem}

When $\mu$ is (the class of) an irreducible representation of $K$, we will write $\mathbf{V}_{G}(\mu)$ for the real-infinitesimal-character irreducible tempered representation of $G$ with lowest $K$-type $\mu$.

\begin{remark} \label{classeEtendue}  We shall need to use Theorem \ref{thmktypes} for groups $G$ that are not quite the group of real points of a connected reductive algebraic group; therefore let us note that Theorem \ref{thmktypes} is also valid when  $G$ is a linear reductive group in Harish-Chandra's class and has all its Cartan subgroups abelian (see \S 0.1 in \cite{Vogan81}). 
 
 \end{remark}
\begin{remark} \label{multiKtype} Still assuming that $G$ is a linear reductive group in Harish-Chandra's class with abelian Cartan subgroups, Vogan proved that if $\pi$ is an irreducible admissible representation, then the lowest $K$-types of $\pi$ all occur with multiplicity one in $\pi|_K$ (see \cite[\S 6.5]{Vogan81}, and for comments, \cite[\S 6]{HermitianFormsSMF}). 
\end{remark}

\begin{remark} \label{SV} The relationship between the highest weight of $\mu$ and the infinitesimal character of $\mathbf{V}_G(\mu)$ is an important ingredient of the results in \cite{Vogan81}. An improved approach, using Carmona's observation that the subtle part of the relationship can be interpreted as a convex cone projection, is explained in \cite{SalamancaRiba_Vogan}, which also contains a summary of the lowest-$K$-type approach to the classification of representations in \cite{Vogan79, Vogan81}. 
\end{remark}


\subsection{Construction of the bijection}\label{subsec:correspondance}

\noindent We first proceed to build tempered representations of $G$ by following as closely as possible the procedure described in \S \ref{subsec:dualG0} to build representations of $G_0$ from Mackey parameters.

Fix a Mackey parameter $(\chi, \mu)$. 

Consider, as we did in \S \ref{subsec:dualG0}, the centralizer $L_\chi$ of $\chi$ in $G$ for the coadjoint action. This is the group of real points of a connected reductive algebraic group. Furthermore, the ``little group'' $K_\chi$ is a maximal compact subgroup of $L_\chi$. Applying Vogan's Theorem \ref{thmktypes}, we can build, out of the irreducible representation $\mu$ of $K_\chi$, a tempered irreducible representation of  $L_\chi$ : the representation $\mathbf{V}_{L_{\chi}}(\mu)$.  %

The centralizer $L_{\chi}$ is the $\theta$-stable Levi factor of a real parabolic subgroup of $G$ (see \cite[Lemma 3.4(4)]{VoganLanglands}). Write $L_{\chi} = M_{\chi} A_{\chi}$ for its Langlands decomposition $-$ a direct product decomposition in which $A_{\chi}$ is abelian and contained in $\exp_{G}(\pe)$, so that $\chi$ defines an abelian character of $L_{\chi}$. 

From our Mackey datum $(\chi, \mu)$, we can thus build as in \S  \ref{subsec:dualG0} a tempered irreducible representation $\sigma = \mathbf{V}_{L_{\chi}}(\mu)\otimes e^{i\chi}$ of the centralizer $L_{\chi}$. 

In order to obtain a tempered irreducible representation of $G$, it is however not reasonable to keep imitating the previous construction and induce without further precaution from $L_{\chi}$ to $G$. Indeed, our centralizer $L_{\chi}=Z_{G}(\chi)$ in $G$ is a poor geometric analogue of the centralizer $L^0_{\chi} = Z_{G_{0}}(\chi)$: these two groups usually do not have the same dimension and the contraction of $L_{\chi}$ with respect to $K_{\chi}$ is usually not isomorphic with $L_{\chi}^0$. The induced representation $\ind_{L_{\chi}}^G(\sigma)$ is also likely to be very reducible.

But we recalled that there exists at least one \emph{parabolic subgroup} $P_{\chi} = L_{\chi}N_{\chi}$ of $G$ whose $\theta$-stable Levi factor is $L_{\chi}$. The parabolic subgroups with $\theta$-stable Levi  factor $L_\chi$ are all contained in a single associate class of parabolic subgroups. 

Now, the subgroup $P_\chi$ is closer to being a good geometric analogue in $G$ of the centralizer $L_{\chi}^0$: these two groups have the same dimension, and whenever $P'_{\chi}$ is a parabolic subgroup of $G$ with $\theta$-stable Levi factor $L_{\chi}$, the contraction of $P'_\chi$ with respect to $K_\chi$ is none other than $L^0_{\chi}$. 

It seems, then, that the representation obtained by extending $\sigma$ to $P_{\chi}$ is the simplest possible analogue of the representation of $L_{\chi}^0$ used in \S \ref{subsec:dualG0}. This makes it natural to form the induced representation
\begin{equation} \label{inductionG} \mathbf{M}(\delta) = \ind_{P_{\chi}}(\sigma) = \ind_{M_{\chi}A_{\chi}N_{\chi}}^{G}(\mathbf{V}_{M_{\chi}}(\mu) \otimes e^{i\chi} \otimes 1).\end{equation}
This is a tempered representation. For the last equality in \eqref{inductionG}, we note that $M_\chi$ satisfies the hypothesis of Remark \ref{classeEtendue} (see \cite[page 142]{Vogan81}), so that Theorem \ref{thmktypes} does apply to $M_\chi$ as well as $L_\chi$; and that $\mathbf{V}_{L_{\chi}}(\mu)$ coincides with the representation of $L_{\chi}$ obtained by extending $\mathbf{V}_{M_{\chi}}(\mu)$ to $L_{\chi}$, because both belong to $\widehat{(L_{\chi})}^{\mathrm{RIC}}$ and have the same lowest $K_{\chi}$-type.\\

We can now state our main result.

\begin{theorem} \label{correspondance} \leavevmode
 \begin{enumerate}[(a)]
 \item For every Mackey parameter $\delta$, the representation $\mathbf{M}(\delta)$ is irreducible and tempered.
 \item Suppose $\delta$ and $\delta'$ are two Mackey parameters. Then the representations $\mathbf{M}(\delta)$ and $\mathbf{M}(\delta')$ are unitarily equivalent if and only if $\delta$ and $\delta'$ are equivalent as Mackey parameters.
 \item By associating to a Mackey parameter $\delta$ the equivalence class of $\mathbf{M}(\delta)$ in $\widehat{G}^{\mathrm{temp}}$, one obtains a bijection between $\mathcal{D}$ and $\widehat{G}^{\mathrm{temp}}$.
  \end{enumerate}\end{theorem}
  
Combined with Mackey's description of $\widehat{G_0}$ in \S \ref{subsec:dualG0}, this yields a bijection
 \begin{equation} \label{corrtemp} \mathcal{M}: \widehat{G_0} \overset{\sim}{\longrightarrow} \widehat{G}^{\mathrm{temp}}.\end{equation} 
We will refer to it as the \emph{Mackey--Higson bijection}.

\begin{proof}[{Proof of Theorem \ref{correspondance}.}]Our result is a simple and straightforward consequence of the classification of tempered irreducible representations of $G$. In other words, it is a simple and straightforward consequence of an immense body of work (the bulk of which is due to Harish-Chandra).

We first remark that given Remark \ref{subsec:data}.(i) and the properties of unitary induction, we may assume without losing generality that $\chi$ lies in $\mathfrak{a}^\star$. 

In that case, $M_{\chi}$ is a reductive subgroup of $G$ that contains the compact subgroup $M = Z_{K}(\mathfrak{a})$ mentioned in \S \ref{subsec:data}, point (ii); in addition, $A_{\chi}$ is contained in $A$, it is related with its own Lie algebra $\mathfrak{a}_{\chi}$ through $A_{\chi} = \exp_{G}(\mathfrak{a}_{\chi})$, and the group $M_{\chi}$ centralizes $A_{\chi}$. We recall that $L_{\chi}$ is generated by  $M$ and the root subgroups for those roots of $(\g_\C, \mathfrak{a}_\C)$ whose scalar product with $\chi$ is zero; in the Langlands decomposition, the subgroup $A_{\chi}$ is by definition the intersection of the kernels of these roots (viewed as abelian characters of $A$).  

In order to choose a precise $N_{\chi}$, fix an ordering on $\mathfrak{a}^\star$: with it comes a notion of positive weight for $(\g, \mathfrak{a}_{\chi})$; the sum of root subspaces corresponding to positive weights yields a subalgebra $\mathfrak{n}_{\chi}$, and we set $N_{\chi} = \exp_{G}(\mathfrak{n}_{\chi})$. But of course another choice for $N_{\chi}$ and $P_\chi$ in  \eqref{inductionG}  would yield an equivalent representation of $G$.

Recall that an element in $\mathfrak{a}^\star$ is called \emph{$\mathfrak{a}_{\chi}$-regular} when its scalar product with every root for the pair $(\g_{\C},\mathfrak{a}_{\chi, \C})$ is nonzero. We now call in a theorem due to Harish-Chandra,  although Harish-Chandra never published it himself; see \cite[Theorem 4.11]{KolkVaradarajan} for semisimple $G$,  \cite[Theorem 14.93]{KnappPrincetonBook} for linear connected reductive $G$, and  \cite[Lemma 3.2(5)]{VoganLanglands}  (with proofs in \cite{SpehVogan}) for the current class of reductive groups:  

 \emph{If $\beta \in \mathfrak{a}^\star$ is $\mathfrak{a_{\chi}}$-regular and if $\eta$ is an irreducible tempered representation of  $M_{\chi}$ with real infinitesimal character, then $\ind_{M_{\chi}A_{\chi}N_{\chi}}^{G}(\eta \otimes e^{i\beta})$ is irreducible.}

To show that $\chi$ is $\mathfrak{a}_{\chi}$-regular, consider a root $\gamma$ of $(\g_{\C},\mathfrak{a}_{\C})$ whose scalar product with $\chi$ is zero. Then $\gamma$ is a root for $(\mathfrak{l}_{\chi, \C}, \mathfrak{a}_{\chi, \C})$ (see the remarks on the structure of $L_\chi$ above), so that $\mathfrak{a}_{\chi}$ is contained in the kernel of $\gamma$ (see the same remarks). Thus, the restriction of $\gamma$ to $\mathfrak{a}_{\chi}$ must be zero, and $\gamma$ cannot be a root of $(\mathfrak{g}_{\C}, \mathfrak{a}_{\chi, \C})$, because the latter are precisely the roots which do not vanish on $\mathfrak{a}_{\chi}$.

Using Harish-Chandra's theorem above, we obtain part (a) in Theorem \ref{correspondance}. Part (b) then follows from the usual criteria for equivalence between parabolically induced representations (see for instance \cite[ Theorem 4.11(i))]{KolkVaradarajan}: {if $M_{1}A_{1}N_{1}$ and $M_{2}A_{2}N_{2}$ are two parabolic subgroups containing a parabolic subgroup $MAN$, if $\beta_{1}$ and $\beta_{2}$ are elements of $\mathfrak{a}_{1}^\star$ and $\mathfrak{a}_{2}^\star$ which are respectively $\mathfrak{a}_{1}$ and $\mathfrak{a}_{2}$-regular, and if $\eta_{1}$ and $\eta_{2}$ are two real-infinitesimal-character irreducible tempered representations of $M_{1}$ and $M_{2}$, then the representations $\mathrm{Ind}_{M_{1}A_{1}N_{1}}^{G}(\eta_{1} \otimes e^{i\beta_{1}})$ and $\mathrm{Ind}_{M_{2}A_{2}N_{2}}^{G}(\eta_{2} \otimes e^{i\beta_{2}}$) are equivalent if and only if there is an element in  $W=W(\g_\C, \mathfrak{a}_\C)$ that sends $\beta_{1}$ to $\beta_{2}$ and $\eta_1$ to $\eta_2$.}

We now turn to (c), which is a simple consequence of the Knapp--Zuckerman classification of irreducible tempered representations \cite{KnappZuckerman, KnappZuckerman2}. Suppose $\pi$ is an arbitrary irreducible tempered representation of $G$; then there exists a parabolic subgroup $P=M_{P}A_{P}N_{P}$ of $G$, a tempered irreducible representation $\sigma$ of $M_{P}$ with real infinitesimal character, and an element $\nu$ of $(\mathfrak{a}_{P})^\star$, so that 
\[ \mathrm{Ind}_{M_{P}A_{P}N_{P}}^G\left( \sigma \otimes e^{i\nu} \right)\]
is irreducible and equivalent with $\pi$. See \cite[ Theorem 3.3]{VoganLanglands}; for proofs, see of course \cite{KnappZuckerman, KnappZuckerman2} in the connected compact-center case, and \cite[\S 5.4]{Shelstad}, for disconnected groups in the present class (see also the comments in \cite{DelormeClozel}).

We need to prove that such a representation is equal to $\mathbf{M}(\delta)$ for some Mackey parameter $\delta$. We first remark that $P$ is associate to a parabolic subgroup contained in $P_{\nu}$: since $\nu$ lies in $(\mathfrak{a}_{P})^\star$, the centralizer  $L_{P}=M_{P}A_{P}$ of $\mathfrak{a}_{P}^\star$ for the coadjoint action is contained in $L_{\nu}$; given the above remarks on the Langlands decomposition and the construction of parabolic subgroups, we deduce that $L_{P}$ is contained in $L_{\nu}$, so that $M_{\nu}$ contains $M_{P}$ and that $A_{\nu}$ is contained in $A_{P}$, and finally that $N_{\nu}$ is conjugate to a subgroup $N'$ contained in $N_{P}$. So $P$ is associate to $M_P A_P N'$, which is contained in $P_{\nu}$.

Write $\tilde{A}$ and $\tilde{N}$ for the subgroups of $G$ whose Lie algebras are the orthogonal complements of $\mathfrak{a}_{\nu}$ and $\mathfrak{n}'$ in $\mathfrak{a}_{P}$ and $\mathfrak{n}_{P}$, respectively. Then $A_{P}=A_{\nu}\tilde{A}$ and $N_{P}=\tilde{N}N'$; from the fact that $\mathfrak{m}_{\nu}$ is orthogonal to $\mathfrak{a}_{\nu} \oplus \mathfrak{n}'$, we deduce that $\tilde{A}$ and $\tilde{N}$ are contained in $M_{\nu}$. Since $A_{P}$ is abelian and $N_{P}$ normalizes $A$, we obain
\[ \mathrm{Ind}^{G}_{M_{P}A_{P}N_{P}} \left( \tau \otimes e^{i\nu} \otimes 1 \right) \simeq \mathrm{Ind}_{(M_{P}\tilde{A} \tilde{N})A_{\nu} N'}^G \left( (\tau \otimes e^{0} )\otimes e^{i\nu} \right). \]
Now,  $\tilde{P} = M_{P}\tilde{A} \tilde{N}$ is a subgroup of $M_{\nu}$ and $M_P\tilde{A}$ is the centralizer of $\tilde{A}$ in $M_{\nu}$, so $\tilde{P}$ is a parabolic subgroup of $M_{\nu}$. Finally, we note that $\tilde{\sigma} = \mathrm{Ind}_{M \tilde{A} \tilde{N}}^{M_{\nu}}(\sigma \otimes e^{0})$  is a tempered representation of $M_{\nu}$, that it has real infinitesimal character (by induction in stages from the definition of real infinitesimal character), and that it is irreducible (otherwise $\pi$ would not be). By induction in stages, we conclude that
 \[\mathrm{Ind}^{G}_{P_{\nu}} \left( \tilde{\sigma} \otimes e^{i\nu} \right) \simeq \mathrm{Ind}_{M_{\nu}A_{\nu}N'}^G \left( \mathrm{Ind}^{M_{\nu}}_{\tilde{P}} \left( \sigma \otimes e^{0} \right) \otimes e^{i\nu} \right):\]
thus $\pi$ is equivalent with $\mathbf{M}(\nu, \mu)$, where $\mu$ is the lowest $K_{\nu}$-type of  the representation $\mathrm{Ind}^{M_{\nu}}_{\tilde{P}}(\sigma \otimes e^{0})$. This concludes the proof.\end{proof}

\begin{remark}  It seems likely that Theorem \ref{correspondance} and its proof are valid in a wider class than used here. Presumably they hold as soon as $G$ satisfies Theorem \ref{thmktypes}, the Harish-Chandra irreducibility criterion above, and the Knapp--Zuckerman classification. Concerning Theorem \ref{thmktypes}, some care is definitely needed:  I am not aware of any reference for nonlinear groups, and there are nonlinear groups in Harish-Chandra's class which do not satisfy the (additional) multiplicity-one property. For $G$ in Harish-Chandra's class, the irreducibility criterion is stated in \cite[Lemma 3.2]{VoganLanglands}, and a version of the Knapp--Zuckerman classification is stated in \cite[Theorem 2.9]{Vogan84Unitarizability}.   \end{remark}


\section{First properties of the bijection} \label{subsec:ktypes}

\noindent In our construction of the Mackey--Higson bijection in \S \ref{subsec:correspondance}, the notion of lowest $K$-type for admissible representations of $G$ was of critical importance. Since $K$ is also a maximal compact subgroup in $G_0$, the notion of lowest $K$-type for admissible representations of $G_0$ can similarly be brought in, and every irreducible admissible representation of $G_0$ has a finite number of lowest~$K$-types. 

We now prove that the Mackey--Higson bijection $\mathcal{M}: \widehat{G_0} \to \widehat{G}^{\mathrm{temp}}$ preserves lowest $K$-types.

\begin{proposition} \label{ktypes} Suppose $\pi_0$ is a unitary irreducible representation of $G_0$. Then  the representation $\pi_0$ of $G_0$ and the representation $\mathcal{M}(\pi_0)$ of $G$ have the same lowest $K$-types. \end{proposition}

\begin{proof} We recall that $\pi_0$ is equivalent with $\mathbf{M}_{0}(\chi, \mu)$ for some Mackey datum $(\chi, \mu)$. From Remark \ref{constructionG0}, we know that as a $K$-module, $\mathbf{M}_{0}(\chi, \mu)$ is isomorphic with $\indk\left(\mu \right)$. So we need to show that the lowest $K$-types of the $G$-representation $\mathbf{M}(\chi, \mu)$ are exactly the irreducible $K$-submodules  of $\ind_{K_{\chi}}^K\left(\mu\right)$ with minimal norm.

Now, the definition of parabolic induction shows that as $K$-modules, $\mathbf{M}(\chi, \mu)$ and  $\indk\left(\mathbf{V}(\mu)|_{K_{\chi}}\right)$ are isomorphic (see the ``compact picture'' description of $\mathbf{M}(\chi, \mu)$ in \cite[\S VII.1]{KnappPrincetonBook}). Of course the latter $K$-module contains $\indk\left(\mu\right)$.

Suppose then that $\alpha$ is a lowest $K$-type in $\indk\left(\mathbf{V}(\mu)|_{K_{\chi}}\right)$, but is not a lowest $K$-type in $\indk\left(\mu\right)$. Then there is an element $\mu_{1}$ of $\widehat{K_{\chi}}$, distinct from $\mu$, such that  $\alpha$ is a lowest $K$-type in $\indk(\mu_{1})$, and since $\alpha$ must appear with multiplicity one in $\indk\left(\mathbf{V}(\mu)|_{K_{\chi}}\right)$ (see Remark \ref{multiKtype}), $\mu_{1}$ must appear with multiplicity one in $\mathbf{V}(\mu)|_{K_{\chi}}$. Because the latter has only one lowest $K_{\chi}$-type, we know that $\norme{\mu_{1}}_{\widehat{K_{\chi}}}$ must be greater than $\norme{\mu}_{\widehat{K_{\chi}}}$.

If $\alpha$ were a lowest $K$-type in $\indk\left(\mathbf{V}(\mu_{1})|_{K_{\chi}}\right)$, the representations $\indk\left(\mathbf{V}(\mu)|_{K_{\chi}}\right)$ and $\indk\left(\mathbf{V}(\mu_{1})|_{K_{\chi}}\right)$ would have a lowest $K$-type in common; but that cannot happen, because of the following reformulation of a result by Vogan: 

\begin{lemma} \label{voganKtypes} Suppose $MAN$ is a cuspidal parabolic subgroup of $G$, and $\mu_{1}$, $\mu_{2}$ are inequivalent irreducible $K \cap M$-modules. Then the representations $\ind_{K \cap M}^K\left(\mathbf{V}(\mu_{1})\right)$ and $\ind_{K \cap M}^K\left(\mathbf{V}(\mu_{2})\right)$ have no lowest $K$-type in common.\end{lemma}

\begin{proof} When $\mathbf{V}(\mu_{1})$ and $\mathbf{V}(\mu_{2})$ are in the discrete series of $M$, this follows from Theorem 3.6 in \cite{VoganClassifying} (see also Theorem 1 in the announcement \cite{Vogan77}, and of course \cite{Vogan79}). In the other cases, this actually follows from the same result, but we need to give some details. 

Assume that both $\mathbf{V}(\mu_{1})$ and $\mathbf{V}(\mu_{2})$ are either in the discrete series or nondegenerate limits of discrete series. Then both $\ind_{MAN}^G\left( \mathbf{V}(\mu_{1}) \right)$ and $\ind_{MAN}^G\left( \mathbf{V}(\mu_{2}) \right)$ are irreducible constituents of some representations induced from discrete series, from larger parabolic subgroups if need be: see \cite[Theorem 8.7]{KnappZuckerman}. If $\ind_{M^*A^*N^*}^G\left(\delta_{1} \right)$ (with $\delta_{1}$ in the discrete series of $M^\star$) contains $\ind_{MAN}^G\left( \mathbf{V}(\mu_{1}) \right)$ as an irreducible constituent, it must contain it with multiplicity one, and the set of lowest $K$-types are $\ind_{M^*A^*N^*}^G\left(\delta_{1} \right)$ is the disjoint union of the sets of lowest $K$-types of its irreducible constituents (which are finite in number): see  \cite[Theorem 15.9]{KnappPrincetonBook}. If $\ind_{MAN}^G\left( \mathbf{V}(\mu_{1}) \right)$ and $\ind_{MAN}^G\left( \mathbf{V}(\mu_{2}) \right)$ are constituents of the same representation induced from discrete series, then the desired conclusion follows; if that is not the case we are now in a position to apply Vogan's result to the two representations induced from discrete series under consideration (Vogan's disjointness-of-$K$-types theorem \emph{is} true of reducible induced-from-discrete-series representations). 

Now, if $\mathbf{V}(\mu_{1})$, in spite of its real infinitesimal character, is neither in the discrete series of $M$ nor a nondegenerate limit of discrete series, then there is a smaller parabolic subgroup $M_{\star} A_{\star} N_{\star}$ and a discrete series or nondegenerate limit of discrete series representation $\epsilon_{1}$ of $M_{\star}$ such that $\ind_{K \cap M}^K(\mathbf{V}(\mu_{1})) = \ind_{K \cap M_{\star}}^K(\epsilon_{1})$. If necessary, we can rewrite $\ind_{K \cap M}^K(\mathbf{V}(\mu_{2}))$ in an analogous way; then we can use Vogan's result again, after some embeddings in reducible representations induced from discrete series as above if necessary. This proves Lemma \ref{voganKtypes}. \end{proof}

Coming back to the proof of Proposition \ref{ktypes}, we now know that if $\alpha$ is a lowest $K$-type in $\indk\left(\mathbf{V}(\mu)|_{K_{\chi}}\right)$, then it cannot be a lowest $K$-type in $\indk\left(\mathbf{V}(\mu_{1})|_{K_{\chi}}\right)$. Let then $\alpha_{1}$ in $\widehat{K}$ be a lowest $K$-type in $\indk\left(\mathbf{V}(\mu_{1})|_{K_{\chi}}\right)$; we note that $\norme{\alpha_1}_{\widehat{K}} < \norme{\alpha}_{\widehat{K}}$. If $\alpha_{1}$ were to appear in $\indk(\mu_{1})$, it would appear in $\indk\left(\mathbf{V}(\mu)|_{K_{\chi}}\right)$, and that cannot be the case because $\alpha$ is already a lowest $K$-type there.

We conclude that there exists $\alpha_{1}$ in $\widehat{K}$ and $\mu_{1}$ in $\widehat{K_{\chi}}$ such that 

\begin{itemize}
\item $\norme{\alpha_{1}}_{\widehat{K}} < \norme{\alpha}_{\widehat{K}}$ and $\norme{\mu_{1}}_{\widehat{K_{\chi}}} > \norme{\mu}_{\widehat{K_{\chi}}}$
\item $\alpha_{1}$ is a lowest $K$-type in $\indk\left(\mathbf{V}(\mu_{1})|_{K_{\chi}}\right)$, but it is not a lowest $K$-type in  $\indk\left(\mu_{1}\right)$. 
\end{itemize}

This seems to trigger an infinite recursion, because the same argument can be used again, beginning with $(\alpha_{1}, \mu_{1})$ instead of $(\alpha, \mu)$; however, there are not infinitely many $K$-types which are strictly lower than $\alpha$. Thus our hypothesis that the lowest $K$-type $\alpha$ in  $\mathbf{M}(\chi, \mu)$ is not a lowest $K$-type in $\ind_{K_{\chi}}^K\left(\mu\right)$ leads to a contradiction. 

To complete the proof, notice (from the compatibility of induction with direct sums) the equality of $K$-modules
\begin{equation} \label{sous_module} \mathbf{M}(\chi, \mu) = \indk\left(\mathbf{V}(\mu)|_{K_{\chi}}\right) = \ind_{K_{\chi}}^K\left(\mu\right) \oplus \tilde{M},\end{equation}
with $\tilde{M}$ induced from a (quite reducible) $K_{\chi}$-module.

We already proved that every lowest $K$-type in  $\mathbf{M}(\chi, \mu)$ must occur in as a lowest $K$-type in $\ind_{K_{\chi}}^K\left(\mu\right)$. Conversely, every lowest $K$-type in $\ind_{K_{\chi}}^K\left(\mu\right)$ \emph{does} occur as a lowest $K$-type in $\mathbf{M}(\chi, \mu)$: from \eqref{sous_module} we know that it must occur there, and from the above argument we already know that it is lower than every $K$-type occuring in $\tilde{M}$, so that it is a lowest $K$-type in $\ind_{K_{\chi}}^K\left(\mu\right)$. This proves Proposition \ref{ktypes}. \end{proof}

We should mention that the above proof establishes, in a somewhat roundabout way, the following property of unitary representations of $G_0$:

\begin{cor} \label{ktypes_G0} {In every unitary irreducible representation of $G_{0}$, each lowest $K$-type occurs with multiplicity one.} \end{cor}

\begin{proof} Suppose $\pi_0$ is a unitary irreducible representation of $G_0$; choose a Mackey parameter $(\chi, \mu)$ for $\pi_0$; recall that $\mathbf{M}_0(\chi, \mu)$ is equivalent, as a $K$-module, with $\ind_{K_{\chi}}^K\left(\mu\right)$. In  \eqref{sous_module} we see that the multiplicity of every lowest $K$-type in $\pi_0$ is lower than its multiplicity in $\mathbf{M}(\chi,\mu)$; every lowest $K$-type occurs at least once in $\pi_0$ and exactly once in $\mathbf{M}(\chi,\mu)$; the corollary follows. \end{proof}

The Mackey--Higson bijection is compatible with other natural features of $\widehat{G}^{\mathrm{temp}}$ and $\widehat{G_0}$. We mention an elementary one, related to the existence in each dual of a natural notion of \emph{renormalization of continuous parameters} for irreducible representations:

\begin{itemize}
\item[$\bullet$] On the $G_0$-side, there is for every $\alpha>0$ a natural bijection 
\begin{equation*} \mathcal{R}^{\alpha}_{G_0}: \widehat{G_0} \to \widehat{G_0}\end{equation*}
obtained from Mackey's description of $ \widehat{G_0} $ in \S \ref{subsec:dualG0} by sending a representation $\pi_0 \simeq \mathbf{M}_0(\chi, \mu)$ of $G_0$ to the equivalence class of $\mathbf{M}_0(\frac{\chi}{\alpha}, \mu)$.

\item[$\bullet$] On the $G$-side, there is also for every $\alpha>0$ a natural bijection 
\begin{equation*} \label{renorm} \mathcal{R}^{\alpha}_G: \widehat{G}^{\mathrm{temp}} \to \widehat{G}^{\mathrm{temp}}\end{equation*}
obtained from the Knapp--Zuckerman classification (see the proof of Theorem \ref{correspondance}) as follows. 

Starting with $\pi$ in $\widehat{G}^{\mathrm{temp}}$, we know that $\pi$ is equivalent with some representation $\mathrm{Ind}_{M_{P}A_{P}N_{P}}^G\left( \sigma \otimes e^{i\nu} \right)$, where $M_{P}A_{P}N_{P}$, $\sigma$ and $\nu$ are as in the proof of Theorem \ref{correspondance}. The tempered representation $\pi^\alpha = \mathrm{Ind}_{M_{P}A_{P}N_{P}}^G\left( \sigma \otimes e^{i\frac{\nu}{\alpha}} \right)$ has the same $R$-group (\cite[\S 10]{KnappZuckerman}) as $\mathrm{Ind}_{M_{P}A_{P}N_{P}}^G\left( \sigma \otimes  e^{i\nu}\right)$, so it is irreducible. We define $\mathcal{R}^{\alpha}_G(\pi)$ as the equivalence class of $\pi^\alpha$.
\end{itemize}

It is immediate from the construction in \S \ref{subsec:correspondance} that the Mackey--Higson bijection is compatible with these renormalization maps:

\begin{proposition} \label{renorm} For every $\alpha >0$, we have $\mathcal{R}^\alpha_G \circ \mathcal{M} =  \mathcal{M} \circ  \mathcal{R}^\alpha_{G_0}$.  \end{proposition}

Together Propositions \ref{ktypes} and \ref{renorm} may shed some light on our construction of the correspondence: any bijection
\[ {\mathcal{B}}:  \widehat{G_0} \to \widehat{G}^{\mathrm{temp}}\]
must, if it is to satisfy Proposition \ref{renorm}, induce a bijection between the fixed-point-sets
\[ \left\{ \pi_0 \in \widehat{G_0}: \forall \alpha>0, \ \mathcal{R}^\alpha_{G_0}(\pi_0)=\pi_0 \right\}  \mathrm{and}    \left\{ \pi \in \widehat{G}^{\mathrm{temp}}: \forall \alpha>0, \ \mathcal{R}^\alpha_{G}(\pi)=\pi \right\}. \]
The first class is the copy of $\widehat{K}$ discussed in Remark \ref{copiekhat}; the second is the class $\widehat{G}^{\mathrm{RIC}}$ discussed in \S \ref{subsec:vogan}. So $ {\mathcal{B}}$ must send any $K$-type $\mu$ to an element of  $\widehat{G}^{\mathrm{RIC}}$; if $ {\mathcal{B}}$ is to satisfy Proposition \ref{ktypes}, it must send $\mu$ to $\mathbf{V}_G(\mu)$, like our correspondence $\mathcal{M}$.

Now, observing the Knapp--Zuckerman classification,  Lemma \ref{voganKtypes} and Proposition \ref{ktypes}, and inserting the above remark, we note that if $\pi_0 \simeq \mathbf{M}(\chi, \mu)$ is a unitary irreducible representation of $G$ with nonzero $\chi$, any representation of $G$ having the same set of lowest $K$-types as $\pi_0$ must read 
\[ \ind_{M_\chi A_\chi N_\chi}\left( \mathbf{V}_{M_\chi}(\mu) \otimes e^{i \nu}\right)\]
for some $\nu$ in $\mathfrak{a}_\chi^\star$. We may view our correspondence $\mathcal{M}$ as that obtained by choosing $\nu=\chi$.


\section{Extension to the admissible duals} \label{subsec:admissible}

\noindent We now come back to the setting of \S \ref{subsec:correspondance} and exend the Mackey--Higson bijection $\mathcal{M}: \widehat{G_0} \to \widehat{G}^{\mathrm{temp}}$ to a natural bijection between the admissible duals $\widehat{G_0}^{\mathrm{adm}}$ and $\widehat{G_{\ }}^{\mathrm{adm}}$.

Let us first describe a way to build, out of an admissible Mackey datum $\delta=(\chi, \mu)$, an admissible representation of $G$ (here $\chi \in \mathfrak{a}_\C^\star$, $\mu \in \widehat{K_\chi}$,  see \S \ref{subsec:admG0}). 

Write $\chi = \alpha+i\beta$, where $\alpha$ and $\beta$ lie in $\mathfrak{a}^\star$. Consider the centralizer $L_\chi$ of $\chi$ for the coadjoint action: $L_\chi$ is the intersection $L_\alpha \cap L_\beta$. Notice that in the notations of the proof of Theorem \ref{correspondance}, $L_\chi$ is the centralizer of $\mathfrak{a}_\alpha + \mathfrak{a}_\beta$; thus $L_\chi$ appears as the $\theta$-stable Levi factor of a real parabolic subgroup $P_\chi=M_\chi A_\chi L_\chi$ of $G$, in which $K_\chi$ is a maximal compact subgroup. In slight contrast with the construction of \S \ref{subsec:correspondance}, we require from the outset that $P_\chi$ be chosen in a specific way: we fix an ordering of $\mathfrak{a}^\star$ making $\alpha$ (the ``real part'' of $\chi$) nonnegative, and choose $P_\chi$ accordingly. See the proof of Theorem \ref{correspondance}.

As before, we can build from $\delta$ the representation $\mathbf{V}_{M_\chi}(\mu) \otimes e^{\chi}$ of $P_\chi$; we define
\begin{equation*} \tilde{\mathbf{M}}^{}(\delta) = \mathrm{Ind}_{P_{\chi}}^G(\sigma) =  \ind_{M_{\chi}A_{\chi}N_{\chi}}^{G}(\mathbf{V}_{M_{\chi}}(\mu) \otimes e^{\chi} \otimes 1),\end{equation*}
an admissible representation of $G$. 

\begin{theorem} \label{corradm} \leavevmode
\begin{enumerate}[(a)]
 \item For every admissible Mackey parameter $\delta$, the representation $\tilde{\mathbf{M}}(\delta)$ has a unique irreducible quotient ${\mathbf{M}^{\mathrm{adm}}}(\delta)$.
 \item Suppose $\delta$ and $\delta'$ are two admissible Mackey parameters. Then ${\mathbf{M}^{\mathrm{adm}}}(\delta)$ and ${\mathbf{M}^{\mathrm{adm}}}(\delta')$ are equivalent if and only if $\delta$ and $\delta'$ are equivalent as admissible Mackey parameters.
 \item By associating to any admissible Mackey parameter $\delta$ the equivalence class of $\mathbf{M}^{\mathrm{adm}}(\delta)$ in $\widehat{G}^{\mathrm{adm}}$, one obtains a bijection between $\mathcal{D}^\mathrm{adm}$ and $\widehat{G}^{\mathrm{adm}}$.
  \end{enumerate}
\end{theorem}

Combined with the description of $\widehat{G_0}^{\mathrm{adm}}$ in \S \ref{subsec:admG0}, this yields a natural bijection 
 \begin{equation} \mathcal{M}^{\mathrm{adm}}: \widehat{G_0}^{\mathrm{adm}} \overset{\sim}{\longrightarrow} \widehat{G}^{\mathrm{adm}}\end{equation}
whose restriction to $\widehat{G_0}$ is the Mackey--Higson bijection \eqref{corrtemp} between tempered representations.\\

\begin{proof}[{Proof of Theorem \ref{corradm}.}] 

We use the Langlands classification to reduce the description of the irreducible admissible representations of $G$ to that of the irreducible tempered representations of reductive subgroups of $G$. 

We will write $(\mathfrak{a}_\alpha)^o$ for the set of linear functionals on $\mathfrak{a}$ whose restriction to $\mathfrak{a}_\alpha$ is zero. Decompose
\[ \chi = \alpha + i (\beta_1 + \beta_2), \quad \text{ where } \quad \alpha \in \mathfrak{a}^\star, \quad \beta_1 \in \mathfrak{a}_{\alpha}^\star \enskip \text{and} \enskip \beta_2 \in (\mathfrak{a}_{\alpha})^o.\]
{Consider the parabolic subgroup $P_\alpha = M_\alpha A_\alpha N_\alpha$ attached to $\alpha$.} We point out that it is possible to view $(\beta_1, \mu)$ as a (tempered) Mackey parameter for the reductive group $M_\alpha$: indeed, $M_\alpha$ admits $K_\alpha$ as a maximal compact subgroup, and if we view $\beta_1$ as a linear functional on $\mathfrak{m}_\alpha \cap \pe$ (where $\mathfrak{m}_\alpha$ is the Lie algebra of $M_\alpha$), then the stabilizer of $\beta_1$ for the action of $K_\alpha$ on $(\mathfrak{m}_\alpha \cap \pe)^\star$ is $K_\alpha \cap K_{\beta_1} = K_\chi$ (since $K_\alpha$ is contained in $K_{\beta_2}$). {Besides, the group $M_\alpha$ satisfies the hypothesis of Remark \ref{classeEtendue}: we can therefore construct, as in \S \ref{subsec:correspondance}, the tempered representation $\mathbf{M}_{M_{\alpha}}(\beta_1, \mu)$.}

We now remark that 
\begin{equation} \label{langl} \tilde{\mathbf{M}}(\chi, \mu) \simeq \mathrm{Ind}_{P_\alpha}^G \left( \mathbf{M}_{M_{\alpha}}(\beta_1, \mu)  \otimes e^{\alpha+i{\beta_2}} \right). \end{equation}
To see this, recall that $\mathbf{M}_{M_{\alpha}}(\beta, \mu)$ is built using the centralizer
\[ Z_{M_\alpha}(\beta_1) = \left\{ g \in M_\alpha, \mathrm{Ad}^\star(g) \beta_1 = \beta_1 \right\} = M_\alpha \cap L_{\beta_1} = M_\chi(A_\chi \cap M_\alpha) \ ;\] thus, it is induced from a parabolic subgroup $P$ of $M_\alpha$, which reads $P=M_\chi(A_\chi \cap M_\alpha)\tilde{N}$ $-$ here $\tilde{N}$ normalizes $Z_{M_\alpha}(\beta_1)$ {and its Lie algebra is a sum of positive root spaces for the given ordering}. The right-hand side of \eqref{langl} can thus be rewritten as 
\[  \ind_{M_\alpha A_\alpha N_\alpha}^{G} \left[ \ind^{M_\alpha}_{M_\chi(A_\chi \cap M_\alpha) \tilde{N}} \left( \mathbf{V}_{M_\chi}(\mu) \otimes e^{i \beta_1} \right) \otimes e^{\alpha + i \beta_2} \right], \]
which, by induction in stages, can in turn be written with only one parabolic induction as
\[  \ind_{M_\chi\left[(A_\chi \cap M_\alpha) (A_\alpha) \right](\tilde{N}N_\alpha)}^{G} \left( \mathbf{V}_{M_\chi}(\mu) \otimes e^{\alpha + i \beta} \right). \]
Since $(A_\chi \cap M_\alpha) (A_\alpha)$ is none other than $A_\chi$, the parabolic subgroup that appears in the final expression is {equal to} $P_\chi$. This proves \eqref{langl}.

From \eqref{langl} we see that $\tilde{\mathbf{M}}(\chi, \mu)$ is induced from a tempered representation; furthermore, in \eqref{langl} we know from the proof of Theorem \ref{correspondance} that $\alpha$ is $\mathfrak{a}_\alpha$-regular. {It also has a positive scalar product with every positive root of $(\mathfrak{g}, \mathfrak{a}_\alpha)$ in the ordering used to define $P_\alpha$}. Therefore $\tilde{\mathbf{M}}(\chi, \mu)$ meets the usual criterion for the existence of a unique irreducible quotient: see \cite[Theorem 8.54]{KnappPrincetonBook}. This proves part (a) in Theorem \ref{corradm}; part (b) then follows from the usual criteria for equivalence between Langlands quotients (implicit in  \cite[Theorem 8.54]{KnappPrincetonBook}).

To prove part (c), recall from the Langlands classification that when $\pi$ is an admissible irreducible representation of $G$, there exists a parabolic subgroup $P=M_P A_P N_P$ of $G$ such that $A_P \subset A$, there exists an $\mathfrak{a}_P$-regular element $\alpha$ in $\mathfrak{a}_P^\star$, and there exists an element $\beta_1$ in $\mathfrak{a}_P^\star$, as well as an irreducible tempered representation $\sigma$ of $M_P$, so that  $\pi$ is equivalent with the unique irreducible quotient of $\mathrm{Ind}_{P}^G \left( \sigma  \otimes e^{\alpha+i{\beta_1}} \right)$. Since $\alpha$ is $\mathfrak{a}_P$-regular, the subgroups $P_\alpha$ and $P$  are associate. Applying Theorem \ref{correspondance} to $M_\alpha$, we know that the tempered representation $\sigma$ of $M_P$ is equivalent with  $\mathbf{M}_{M_\alpha}(\beta_2, \mu)$ $-$ for some linear functional  $\beta_2$ on $\mathfrak{a} \cap (\mathfrak{m}_\alpha \cap \pe)$ and some irreducible representation $\mu$ of $K_\alpha \cap Z_{M_\alpha}(\beta_2) = K_{\alpha+i\beta_1+i\beta_2}$. (The group $M_\alpha$ does not necessarily satisfy the ``real points of a connected complex algebraic group'' hypothesis of \S \ref{subsec:correspondance}, but the centralizer $L_\alpha$ does, and it is the direct product between $M_\alpha$ and the abelian vector group $A_\alpha$; thus Theorem \ref{correspondance}  does~hold~for~$M_\alpha$). We note that $\beta_2 \in (\mathfrak{a}_\alpha)^o$, and conclude that $\pi$ is equivalent with a representation of the kind that appears in the right-hand side of \eqref{langl}. This proves (c) and Theorem \ref{corradm}.  \end{proof}


\providecommand{\bysame}{\leavevmode\hbox to3em{\hrulefill}\thinspace}
\providecommand{\MR}{\relax\ifhmode\unskip\space\fi \textcolor{blue}{MR}}
\providecommand{\MRhref}[1]{%
  \href{http://www.ams.org/mathscinet-getitem?mr=#1}{#1}
}
\providecommand{\href}[2]{#2}

\end{document}